\documentclass[letterpaper, 10 pt, conference]{ieeeconf}  % Comment this line out if you need a4paper
\IEEEoverridecommandlockouts                              % This command is only needed if 
\usepackage[utf8]{inputenc} % allow utf-8 input
\usepackage[T1]{fontenc}    % use 8-bit T1 fonts
\usepackage{hyperref}       % hyperlinks
\usepackage{url}            % simple URL typesetting
\usepackage{booktabs}       % professional-quality tables
\usepackage{amsfonts}       % blackboard math symbols
\usepackage{nicefrac}       % compact symbols for 1/2, etc.
\usepackage{microtype}      % microtypography
\usepackage{graphicx}
\usepackage{amsmath}
\usepackage{amssymb}
\usepackage{colortbl}
\usepackage{mathtools}
\usepackage{flushend}
\usepackage{cite}
\usepackage{siunitx}
\usepackage{multirow}

\definecolor{Gray}{gray}{0.9}

\title{Constrained Block Nonlinear Neural Dynamical Models}

\author{Elliott Skomski, Soumya Vasisht, Colby Wight, Aaron Tuor, J\'an Drgo\v na, Draguna Vrabie \\
    Pacific Northwest National Laboratory,
	Richland, Washington, USA,\\
	\{elliott.skomski, soumya.vasisht, colby.wight, aaron.tuor, jan.drgona, draguna.vrabie\}@pnnl.gov
}
\begin{document}

\maketitle

\begin{abstract}
Neural network modules conditioned by known priors can be effectively trained and combined to represent systems with nonlinear dynamics. This work explores a novel formulation for data-efficient learning of deep control-oriented nonlinear dynamical models by embedding local model structure and constraints. The proposed method consists of neural network blocks that represent input, state, and output dynamics with constraints placed on the network weights and system variables. 
For handling partially observable dynamical systems, we utilize a state observer neural network to estimate the states of the system's latent dynamics. 
We evaluate the performance of the proposed architecture and training methods on system identification tasks for three nonlinear systems: a continuous stirred tank reactor, a two tank interacting system, and an aerodynamics body. Models optimized with a few thousand system state observations accurately represent system dynamics in open loop simulation over thousands of time steps from a single set of initial conditions. Experimental results demonstrate an order of magnitude reduction in open-loop simulation mean squared error for our constrained, block-structured neural models when compared to traditional unstructured and unconstrained neural network models. 

\end{abstract}

\section{Introduction}
% Soumya  \\
Most systems are inherently nonlinear in nature, which makes the study of nonlinear models of systems of utmost importance to engineers, physicists, and mathematicians. In contrast with linear systems, where observations are proportional to changes in the input and state variables, nonlinear systems appear unpredictable, counter-intuitive, or even chaotic \cite{khalil_2002}. Identifying mathematical models for such complex dynamical systems from input-output data is central to most engineering systems for modeling, prediction and control. 

Most works in deep learning-based system identification train large, complex, monolithic networks with strong regularization to avoid overfitting. However, greater data efficiency and stability in optimization of more precise and less opaque systems models can be achieved by embedding prior knowledge and intuition of general system characteristics such as structure and dissipation \cite{ioann_thesis}. In addition, modern neural network approaches don't typically provide constraints enforcement within physically realizable and operationally safe bounds, which can lead to unexpected model behavior in unobserved regions of the system state space. To address these shortcomings, we propose a ``domain intuitive'' framework for constructing constrained block nonlinear neural network dynamics models which incorporate general knowledge and intuition about system behavior and can train efficiently from small data sets. 
% In an effort to drift away from these obscure, unstructured models, embedding prior knowledge about certain system characteristics has been shown to improve model learning and modeling accuracy \cite{ioann_thesis}.

We present a generalized family of neural state space model architectures in the form of constrained block-structured nonlinear neural state space models. Our approach combines the benefits of classical system identification techniques with the advances in constrained deep learning to learn deep structured models. In particular, we build modular representations of nonlinear systems with identifiable neural modules that may be constrained and influenced by structural priors. This allows for inclusion of stability constraints through eigenvalue regularizations of the network weights, state and input constraints, and other boundary condition constraints.

We investigate the benefits of the proposed method on three nonlinear system identification tasks,  and analyze, through an ablation study, the value of adding structure and constraints on the modeling outcome. Compared to their unstructured unconstrained counterparts, the presented constrained block-structured models show 68-91\% reduction in open-loop MSE from the ground truth dynamics. Our experimental results %\footnote{For reproducing the research, we make our open-source code freely available (url to be included).} 
empirically demonstrate the potential to learn physically representative neural network dynamics models from limited recorded systems measurements, and further underscore the significance of including structure and constraints in the identification process.

\section{Background}
 
 Various methods exist for identifying nonlinear systems differing in their purpose and application \cite{schoukens2019}. When a white-box approach of mathematical modeling using first principles is not feasible, data-driven gray and black-box methods have been proposed including Volterra series, neuro-fuzzy models, block-oriented models, and NARMAX models \cite{bai_introduction_2010}. These methods have been successful in many areas \cite{brockett_volterra_1976, oh2011application, chiras2001}. 

In the current work we focus on block-oriented models which are composed of two basic building blocks: linear dynamic models and static nonlinearities \cite{schoukens2019}. Such block-structured nonlinear models, represent a nonlinear system as a collection of static or dynamic components, either linear or simpler nonlinear sub-systems \cite{bai_introduction_2010}. This enables a more accurate description of the process behavior and reduces the hard problem of identifying high-order nonlinear systems to the identification of lower order subsystems and their interactions. Two such simple models are the Hammerstein and the Wiener models. 
Hammerstein systems are useful to describe linear dynamic processes driven by a nonlinear actuator. Wiener systems describe linear dynamic processes equipped with a nonlinear sensor. Despite their simplicity, such models have been successfully used in many engineering fields (signal processing, identification of biological systems, modeling of distillation columns, modeling of hydraulic actuators, etc.) \cite{kibangou2006, gilabert2005, Bai2003, dempsy2004, kwak1998}, providing reasonable approximations and insight into the system structure \cite{pearson_gray-box_2000}.

Neural networks have a long  and productive history in data-driven physical modeling. Neural networks were used as early as 1990 to model discrete-time nonlinear dynamics \cite{chen1990}.
Recurrent neural networks (RNN) and long-short term memory (LSTM) networks are dynamics-inspired architectures which have been studied heavily in the context of machine learning applications
\cite{chow1998recurrent, al2008nonlinear, yu2004nonlinear, de2007nonlinear}. 
 In recent years, RNNs and LSTMs are becoming increasingly popular for their potential in tackling nonlinear problems in black-box system identification~\cite{LSTM_SysID2017}. For instance, an LSTM-based architecture in \cite{woo_dynamic_2018} was used to model nonlinear effects in unmanned surface vehicle dynamics. 
 Several extensions to straightforward applications of recurrent architectures have since been proposed.
 For instance, several variations of Neural Kalman filters \cite{wilson_neural_2009, coskun2017long, krishnan2015deep} have been proposed for control oriented systems modeling, and neural Koopman operators~\cite{yeung2019learning} have been presented. In addition, several recent advances have extended neural networks to model nonlinear continuous time dynamics \cite{chen2018neural, champion2019data, Raissi_PINN_p1_2017, raissi_multistep_2018}. 
 
Some recent works interpret deep neural networks through the optics of differential equations~\cite{HeZRS15, NIPS2018_7892_NeuralODEs, DeepXDE2019} creating stronger theoretical links between dynamical systems and neural networks.
Others show that structural priors via factored parametrizations of linear maps can be used, for instance to constrain the eigenvalues of the layer Jacobians~\cite{IMEXnet2019}.
Additionally, it has been demonstrated that regularizing the eigenvalues of successive linear transformations alleviates 
the vanishing and exploding gradient problem, allowing to train very deep recurrent architectures~\cite{HaberR17}.

Neural state space models (SSM)~\cite{krishnan2016structured,LatentDynamics2018,MastiCDC2018,NIPS2018_8004} represent structural modifications of vanilla RNNs. For example, authors in \cite{ogunmolu_nonlinear_2016} compared different neural network architectures including strictly feed-forward and recurrent Hammerstein models for identifying nonlinear dynamical systems. Our present work is a generalization of several prior works on structured neural SSMs, such as Hammerstein~\cite{OgunmoluGJG16}, or Hammerstein-Wiener models~\cite{HW_RNN2008}. 

\section{Methods} \label{sec: models}

Consider the non-autonomous partially observable nonlinear dynamical system
\begin{subequations}
\label{eq:nonautosystem}
\begin{align}
  & \frac{d}{dt}\mathbf{x}(t) = \mathbf{f}(\mathbf{x}(t),\mathbf{u}(t)), \\ 
  & \mathbf{y}(t) = \mathbf{g}(\mathbf{x}(t)), \\ 
  & \mathbf{x}(t_0) = \mathbf{x}_0, \ \ \ \mathbf{u}(t_0) = \mathbf{u}_0
\end{align}
\end{subequations}
where $\mathbf{y}(t) \in \mathbb{R}^{n_y}$ are  measured outputs, and $\mathbf{u}(t) \in \mathbb{R}^{n_u}$ represents controllable inputs or measured disturbances affecting the system dynamics. 

\subsection{Neural network dynamical models}
\subsubsection{Unstructured nonlinear dynamical models}
In the absence of prior knowledge of the full dynamics of the system and the interactions between its components, we use the following generic form of a neural state space model:
% TODO(lltt): i think these are essentially the same as the previous eqns?
\begin{subequations}
    \label{eq:blackboxssm}
\begin{align}
    \mathbf{x}_{t+1} &= \mathbf{f}_{xu}([\mathbf{x}_t; \mathbf{u}_t]) \\
    \mathbf{\hat{y}}_{t+1} &= \mathbf{f}_{y}(\mathbf{x}_{t+1})
\end{align}
\end{subequations}
where $\mathbf{f}_{xu}$ is a neural network which models the interaction between state and input/disturbance vectors, and $\mathbf{f}_y$ is a linear map. Although this formulation provides a good substrate for learning system dynamics with unknown properties, it lacks structure for incorporating functional priors for specific system components. 

\subsubsection{Block-nonlinear dynamical models}
\begin{table}[t]
    \centering
    \caption{Linear components for block model classes.}
    \begin{tabular}{lccc}
        \toprule
      {}  & \multicolumn{3}{c}{Linear Map}\\
        Model class & $\mathbf{f}_x$ & $\mathbf{f}_u$ & $\mathbf{f}_y$ \\
        \midrule
        Block-nonlinear & N & N & Y \\
        Hammerstein-Wiener &  Y & N &  N\\
        Hammerstein & Y & N & Y\\
        Wiener & Y & Y & N\\
        Linear & Y & Y & Y\\
        \bottomrule
    \end{tabular}
    \label{tab:block}
    \vspace{-10pt}
\end{table}

Given prior assumptions of the interactions between the components of a block-structured nonlinear model, we can instead use a formulation with a structured representation of the state and input/disturbance dynamics of a system:
\begin{subequations}
    \label{eq:blocknlinssm}
\begin{align}
    \mathbf{x}_{t+1} &= \mathbf{f}_x(\mathbf{x}_t) + \mathbf{f}_u(\mathbf{u}_t) \\
    \mathbf{\hat{y}}_{t+1} &= \mathbf{f}_y(\mathbf{x}_{t+1})
\end{align}
\end{subequations}
where $\mathbf{f}_x$ and $\mathbf{f}_u$ are nonlinear maps, and $\mathbf{f}_y$ is a linear map. This formulation allows each component of the state dynamics to assume separate priors, and for interactions between components to be specified directly. 

By changing the linearity of model components $\mathbf{f}_x$, $\mathbf{f}_u$, or $\mathbf{f}_y$ the formulation can be adjusted to follow the structure of other block-oriented state space models, which may be more appropriate to model particular systems. Table \ref{tab:block} lists the possible model class configurations and the linearity of their respective components.

\subsubsection{State observer model}
To model partially observable dynamical systems, we assume that the states $\mathbf{x}_t$ of our neural state space models represent latent dynamics. Hence, we include an additional neural network representing a state observer:
\begin{equation}
    \label{eq:obsv}
    \mathbf{x}_{0} =  \mathbf{f}_o([\mathbf{y}_{1-N_p}; \ldots; \mathbf{y}_{0}])
\end{equation}
where $\mathbf{f}_o$ is a nonlinear map, and $N_p$ is a lookback horizon capturing potential time lag of observed system outputs ${\mathbf{y}}_{t}$.

\subsubsection{Overall Model Architecture}
\begin{figure}[!htbp]
    \centering
  \includegraphics[width=\columnwidth]{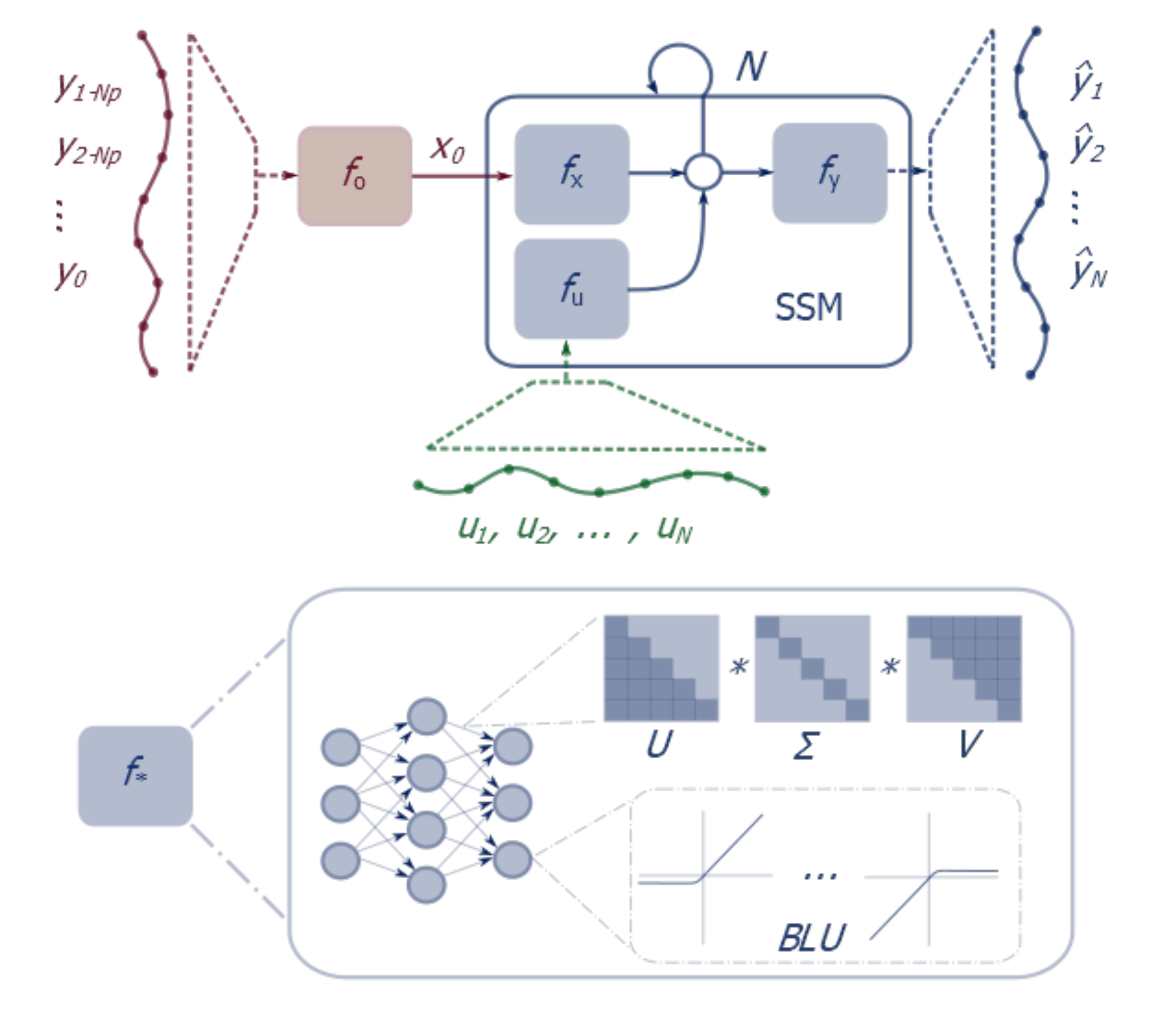}
\caption{Block nonlinear neural state space  model architecture for partially observable systems with $N_p$-step past and $N$-step prediction horizons.
Individual blocks can be parametrized by different structured linear maps with learnable activation functions.}
    \label{fig:SSM}
    \vspace{-10pt}
\end{figure}

Fig.~\ref{fig:SSM} shows the overall architecture of the block structured neural state  model~\eqref{eq:blocknlinssm} with a state estimator model~\eqref{eq:obsv}. Please note that different neural network component blocks can be parametrized by structured linear maps and learnable activation functions explained in the following sections.

%%%%%%%%%%%%%%%%%%%%%%%%%%%%%%%%%%%%%%%%%%%%%%%%%%%%%%%%%%%%%%%%%%%%%%%%%%%%%%%%%%%%%%%%%%%%%%%%%%%%%%%%%%%%%
%%%%%%%%%%%%%%%%%%%%%%%%%%%%%%
%%%%%%%%%%%%%%%%%%%%%%%%%%%%%%%%%%%%%%%%%%%%%%%%%%%%%%%%%%%%%%%%%%%%%%%%%%%%%%%%%%%%%%%%%%%%%%%%%%%%%%%%%%%%
\subsection{Optimization Problem}
We assume access to a limited set of system measurements in the form of tuples, each of which corresponds to the input-output pairs along sampled trajectories with temporal gap $\Delta$. That is, we form a dataset
\begin{equation}
    \label{eq:dataset}
    S = \{(\mathbf{u}^{(i)}_t, \mathbf{y}^{(i)}_t), (\mathbf{u}^{(i)}_{t+\Delta}, \mathbf{y}^{(i)}_{t+\Delta}),  \ldots, (\mathbf{u}^{(i)}_{t+N\Delta}, \mathbf{y}^{(i)}_{t+N\Delta})\},
\end{equation}
where $i = 1,2,\ldots, n$ represents up to $n$ different batches of input-output trajectories with $N$-step time horizon length.  %$n$ gives the total number of batches.

For each batch, $i$, of observations, the neural state estimator generates an initial state $\mathbf{x}_0$ and the neural state space model, in open-loop, recurrently generates $N$ predictions, $\hat{\mathbf{y}}^{(i)}_{1}, ..., \hat{\mathbf{y}}^{(i)}_{N}$. 
The neural network state space model is then fit to the unknown dynamics~\eqref{eq:nonautosystem} by minimizing the mean squared error, $\mathcal{L}_y$, between predicted values and the ground truth measurements for the $N$-step prediction horizon:
\begin{equation}
\label{eq:loss:y}
    \mathcal{L}_y = \frac{1}{n \cdot N \cdot n_y} \sum_{i=1}^n \sum_{t=1}^{N} \sum_{j=1}^{n_y} (\hat{\mathbf{y}}^{(i)}_{t, j} - \mathbf{y}^{(i)}_{t, j})^2
\end{equation}
where $n_y$ is the dimension of the $\mathbf{y}_t$ vectors of observations and predictions.

%%%%%%%%%%%%%%%%%%%%%%%%%%%%%%%%%%%%%%%%%%%%%%%%%%%%%%%%%%%%%%%%%%%%%%%%%%%%%%%%%%%%%%%%%%%%%%%%%%%%%%%%%%%%%
%%%%%%%%%%%%%%%%%%%%%%%%%%%%%%
%%%%%%%%%%%%%%%%%%%%%%%%%%%%%%%%%%%%%%%%%%%%%%%%%%%%%%%%%%%%%%%%%%%%%%%%%%%%%%%%%%%%%%%%%%%%%%%%%%%%%%%%%%%%
\subsection{Neural network block components}
We consider three different classes of neural networks as options for the state observer and nonlinear block components of the state space models. 
\subsubsection{Multi-layer Perceptron}
The MLP is a universal function approximator which composes successive affine linear transformations with a nonlinear activation function $\mathbf{g}$. 
\begin{subequations}
\label{eq:mlp}
\begin{align}
    \mathbf{f} = \mathbf{f}^{(L)} \circ \mathbf{f}^{(L-1)} \circ \hdots  \circ \mathbf{f}^{(1)}\\
    \mathbf{f}^{(k)}(\mathbf{x}) = \mathbf{g}(\mathbf{W}^{(k)}\mathbf{x} + \mathbf{b}^{(k)})
\end{align}
\end{subequations}
where $\mathbf{W}$ and $\mathbf{b}$ are weight and bias terms, respectively. 
The number of hidden layers of the MLP is the number of nonlinear compositions, and the number of nodes for each layer is the dimension of the projected space for the linear map for that layer. Deeper (more layers) and wider (more nodes) networks allow for the approximation of arbitrary multivariable vector-valued functions. 
\subsubsection{Residual MLP}
The residual variant of the MLP adds so called shortcut connnections from the output of the previous layer:
\begin{equation}
    \mathbf{f}^{(k)}(\mathbf{x}) = \mathbf{g}\bigl(\mathbf{W}^{(k)}\mathbf{x} + \mathbf{b}^{(k)} + \mathbf{f}^{(k-1)}(\mathbf{x})\bigr)
\end{equation}
Residual MLPs mitigate vanishing gradients during backpropagation, allowing for effective optimization of deeper neural networks \cite{HeZRS15}.

\subsubsection{Recurrent Neural Networks}
To model sequential dependencies across time, we make use of RNNs which add a recurrent term to the standard MLP as an additional affine linear transformation of the previous transformed time step:
\begin{equation}
    \mathbf{f}^{(k)}(\mathbf{x}_t) = 
        \mathbf{g}\bigl(\mathbf{W}^{(k)}\mathbf{x}_t
         + \mathbf{W}_r^{(k)} \mathbf{f}^{(k)}(\mathbf{x}_{t-1}) + \mathbf{b}^{(k)}\bigr)
\end{equation}
where $\mathbf{W}_r$ are the weights of the recurrent term.

\subsubsection{Activation functions}
We consider two prospective nonlinear activation functions for the aforementioned neural network block functions. The Gaussian Error Linear Unit (GELU)~\cite{HendrycksG16} weights activations by their value
rather than gating them by their sign and has demonstrated superior performance to other activations on 
some tasks. Bendable Linear Units (BLU)~\cite{godfrey2019evaluation} learns an activation from a family of smooth approximations to two piece linear functions such as 
encompassed by Parametric ReLU~\cite{he2015delving} but which are $C^{\infty}$ continuous, yielding a smoother error surface. 

\subsection{Linear maps}

In addition to standard neural network weights, we use structured matrices and regularization strategies to constrain the models' constituent linear maps to satisfy various structural and spectral properties. Constrained linear maps provide a powerful means of embedding prior information into a dynamical model, and can improve the stability of both the underlying optimization problem and the learned dynamics model.

\subsubsection{Perron-Frobenius Parametrization}
We impose constraints on the eigenvalues of a linear map using insight from the Perron-Frobenius theorem: the row-wise minimum and maximum of any positive square matrix defines its dominant eigenvalue's lower and upper bound, respectively. To apply this, we introduce a second weight matrix $\mathbf{M} \in \mathbb{R}^{n_x \times n_x}$ whose values are clamped between lower and upper eigenvalue bounds $\lambda_{\text{min}}$ and $\lambda_{\text{max}}$, then elementwise multiply it with a weight matrix $\mathbf{W} \in \mathbb{R}^{n_x \times n_x}$ with row-wise softmax applied such that each of its rows sum to 1:
\begin{subequations}
\begin{align}
    \mathbf{\Tilde{M}} &= \lambda_{\text{max}} - (\lambda_{\text{max}} - \lambda_{\text{min}}) \cdot \sigma(\mathbf{M})\\
    \mathbf{\Tilde{W}}_{i,j} &= \frac{\text{exp}(\mathbf{W}_{i,j})}{\sum_{k=0}^{n_x} \text{exp}(\mathbf{W}_{i,k})} \cdot \mathbf{\Tilde{M}}_{i,j}
\end{align}
\end{subequations}
where $\sigma(x) = 1 / 1 - \text{exp}(-x)$ is the logistic sigmoid function. 

\subsubsection{Soft SVD Parametrization}
This technique parametrizes a linear map as a factorization via singular value decomposition. The map is defined as two unitary matrices $\mathbf{U}$ and $\mathbf{V}$ initialized as orthogonal matrices, and singular values $\Sigma$ initialized randomly. During optimization, $\mathbf{U}$ and $\mathbf{V}$ are constrained via regularization terms such that they remain orthogonal:
\begin{subequations}
\begin{align}
    \mathcal{L}_{\mathbf{U}} &= || \mathbf{I} - \mathbf{UU}^T ||_2 + || \mathbf{I} - \mathbf{U}^T\mathbf{U} ||_2 \\
    \mathcal{L}_{\mathbf{V}} &= || \mathbf{I} - \mathbf{VV}^T ||_2 + || \mathbf{I} - \mathbf{V}^T\mathbf{V} ||_2 \\
    \mathcal{L}_{\text{reg}} &= \mathcal{L}_{\mathbf{U}} + \mathcal{L}_{\mathbf{V}}
\end{align}
\end{subequations}
Similar to Perron-Frobenius regularization, this regularization also enforces boundary constraints on the singular values $\Sigma$, and its eigenvalues by extension. This is again achieved by clamping and scaling $\Sigma$:
\begin{subequations}
\begin{align}
    \mathbf{\Tilde{\Sigma}} &= \text{diag}(\lambda_{\text{max}} - (\lambda_{\text{max}} - \lambda_{\text{min}}) \cdot \sigma(\Sigma))\\
    \mathbf{\Tilde{W}} &= \mathbf{U\Tilde{\Sigma}V}
\end{align}
\end{subequations}
where $\lambda_{\text{min}}$ and $\lambda_{\text{max}}$ are the lower and upper singular value bounds, respectively.

\subsubsection{Spectral Parametrization}
Similar to Soft SVD, this method proposed by \cite{zhang2018stabilizing} uses an SVD factorization to parametrize a linear map to enforce spectral constraints. This is achieved structurally using Householder reflectors to represent unitary matrices $\mathbf{U}$ and $\mathbf{V}$. The singular values $\Sigma$ are also constrained to lie between $\lambda_{\text{min}}$ and $\lambda_{\text{max}}$ using a sigmoid clamping similar to Soft SVD. In this form, no regularization term is needed to maintain orthogonal structure or spectral properties. 
% For brevity, we refer readers to the original work for further details \cite{zhang2018stabilizing}. % TODO: should we defer to the original paper or go into detail here? it seems like there's a lot to go through.

\subsection{Constraints}
\subsubsection{Output Constraints}
We apply inequality constraints on predictions during training in order to promote the boundedness and convergence of our dynamical models. We define output lower and upper bounds $\mathbf{\underline{y}}$ and $\mathbf{\overline{y}}$, respectively, and derive corresponding slack variables $\mathbf{s^{\underline{y}}}$ and $\mathbf{s^{\overline{y}}}$ to include in the objective function as an additional term $\mathcal{L}^{\text{con}}_y$:
\begin{subequations} \label{eqns:state_obs_const}
\begin{align}
        \mathbf{s^{\underline{y}}} &= \text{max}(0, -\mathbf{\hat{y}} + \mathbf{\underline{y}})\\
        \mathbf{s^{\overline{y}}} &= \text{max}(0, \mathbf{\hat{y}} - \mathbf{\overline{y}})\\
        \mathcal{L}^{\text{con}}_y &= \frac{1}{n_y}\sum_{i=1}^{n_y} (\mathbf{s}^{\mathbf{\underline{y}}}_i + \mathbf{s}^{\mathbf{\overline{y}}}_i)
\end{align}
\end{subequations}
Bounding system trajectories is a straightforward way 
of including prior knowledge by delineating a physically meaningful phase space of the learned dynamics.

\subsubsection{Input Influence Constraints}
For block nonlinear models, we apply additional constraints to control the influence of learned input dynamics $\mathbf{f}_u$ on predictions. Like the state observation constraints previously described, these are inequality constraints and are derived for $\mathbf{f}_u$ in the same manner as in Equation \ref{eqns:state_obs_const} to create another term $\mathcal{L}^{\text{con}}_{\mathbf{f}_u}$. 

\subsubsection{Predicted State Smoothing}
To promote continuous trajectories of our dynamics models, we optionally apply a state smoothing loss which minimizes the mean squared error between successive predicted states:
\begin{equation}
    \mathcal{L}_{dx} = 
        \frac{1}{(N-1)n_x} \sum_{t=1}^{N-1} 
         \sum_{i=1}^{n_x} ({\mathbf{x}}^{(i)}_{t} - {\mathbf{x}}^{(i)}_{t+1})^2
\end{equation}

\subsubsection{Multi-objective loss function}
We include constraints penalties as additional terms to the optimization objective \ref{eq:loss:y}, and further define coefficients, $Q_{*}$ as hyperparameters to scale each term in the multi-objective loss function:
%We define coefficients as hyperparameters to scale each term in the objective function: $Q_y$ for predicted state error, $Q_{\text{reg}}$ for regularization terms, $Q_{dx}$ for state smoothing, $Q^{\text{con}}_x$ for state observation constraints, and $Q^{\text{con}}_{\boldsymbol{f}_u}$ for input/disturbance influence constraints.
%We define coefficients as hyperparameters to scale each term in the objective function: $Q_y$ for predicted state error, $Q_{\text{reg}}$ for regularization terms, $Q_{dx}$ for state smoothing, $Q^{\text{con}}_x$ for state observation constraints, and $Q^{\text{con}}_{\boldsymbol{f}_u}$ for input/disturbance influence constraints.
\begin{equation}
\label{eq:loss:multi}
    \mathcal{L} = Q_y\mathcal{L}_y + Q_{\text{reg}}\mathcal{L}_{\text{reg}} + Q_{dx}\mathcal{L}_{dx} + Q^{\text{con}}_y\mathcal{L}^{\text{con}}_{y} + Q^{\text{con}}_{\mathbf{f}_u}\mathcal{L}^{\text{con}}_{\mathbf{f}_u}
\end{equation}
Note that in this work, $\mathcal{L}_{\text{reg}}$ is only applied if Soft SVD linear map regularization is used. To minimize this objective, we use gradient descent via the AdamW optimizer \cite{loshchilov2017decoupled}.

\section{Numerical Case Studies}
In this section we outline the experimental setup, procedure, and resulting outcomes for three numerical case studies which apply our constrained neural block nonlinear modeling approach to three nonlinear system identification tasks. 
\subsection{Experimental Setup}
We consider three systems with varying degrees and types of nonlinearities in states and inputs to illustrate the versatility of our neural block-structured framework. 

% L_MSE expression
% can put system equations here just in case, if we hit page limit move to appendix
In chemical engineering, the \textit{Continuous Stirred Tank Reactor (CSTR)} is a common mathematical model for a chemical reactor. It is equipped with a mixing device to provide efficient mixing of materials. It is a non-dissipative system that exhibits nonlinear behavior during exothermic reactions with unstable, oscillatory modes that makes the identification process non-trivial.
The model is described by the following set of ordinary differential equations:
\begin{subequations}
\label{eq:CSTR}
\begin{align}
    & r = k_0  e^{-\frac{E}{R\mathbf{x}_2}}  \mathbf{x}_1 \\
    & \dot{\mathbf{x}}_1 = \frac{q}{V}  (C_{af} - \mathbf{x}_1) - r \\
    &  \dot{\mathbf{x}}_2 = \frac{q}{V}  (T_f - \mathbf{x}_2)  
   + \frac{H}{\rho  c_p}  rA + \frac{A}{V \rho c_p  (\mathbf{u} - \mathbf{x}_2)} 
\end{align}
\end{subequations}
where the measured system states $\mathbf{x}_1$, and $\mathbf{x}_2$ are the concentration of the product and temperature, respectively. The system has one control variable $\mathbf{u}$ that represents the temperature of cooling jacket.
We assume that the structure of the differential equations and the parameters $\{k_0, E, R, q, V, H, \rho, c_p, A, T_f, C_{af}\}$ are unknown. 
We generate a training dataset with an emulator based on the system equations. Random step inputs are used to excite the system to capture the system's stability region and transient response. We generate the dataset by simulating the system for 10,000 time steps using the Scipy ODEInt solver.
% \footnote{\url{https://docs.scipy.org/doc/scipy/reference/generated/scipy.integrate.ode.html}}.
% Discuss with Elliott about what model structures were used

The \textit{Two Tank system} we consider consists of two water tanks that are connected by a valve. Liquid inflow to the first tank is governed by a pump. The valve opening controls the flow between the tanks and can either be fully open or fully closed. The system can be described by the ordinary differential equations:
\begin{subequations}
\label{eq:2tank}
\begin{align}
    & \dot{\mathbf{x}}_1 = 
    \begin{cases}
   (1 - \mathbf{u}_1) c_1 \mathbf{u}_2  - c_2 \sqrt{\mathbf{x}_1},  & \text{if } \mathbf{x}_1 \leq 1\\ 
    0 &  \text{otherwise}
    \end{cases} \\
    &  \dot{\mathbf{x}}_2 = 
            \begin{cases}
     c_1  \mathbf{u}_1  \mathbf{u}_2 + c_2  \sqrt{\mathbf{x}_1} - c_2 \sqrt{\mathbf{x}_2} ,  & \text{if } \mathbf{x}_2 \leq 1\\ 
    0 &  \text{otherwise}
    \end{cases} 
\end{align}
\end{subequations}
where the measured system states $\mathbf{x}_1$, and $\mathbf{x}_2$ denote the liquid levels in first and second tank, respectively. The system has two control variables $\mathbf{u}_1$ and $\mathbf{u}_2$ representing the pump speed and valve opening, respectively. Similar to the CSTR system we generate a dataset of 10,000 time steps via simulation using the Scipy ODEInt solver.

The \textit{Aerodynamic Body} exhibits oscillatory behavior when subject to symmetric disturbances. We use a model that predicts the acceleration and velocity of the body using the measurements of its velocities (translational and angular) and various angles related to its control surfaces. It has a large number of parameters, including 4 aerodynamic force coefficients, 11 aerodynamic momentum coefficients, 3 moments of inertia factors, and 4 other constant parameters. The system has 10 inputs, including the various angles of the control surfaces, angle of attack and sideslip, and the measured angular velocities about the three inertial axes. Its 5 outputs measure the angular velocities and accelerations about the lateral y and z direction. For further details about the model, we refer the reader to the dataset source~\cite{aero}. The dataset contains 501 input-output measurements sampled uniformly at 0.02 second intervals.

For each system, we group the time series of inputs and observed measurements into evenly split contiguous training (first \nicefrac{1}{3} of the data points), validation (next \nicefrac{1}{3} of the data points), and test sets (final \nicefrac{1}{3} of the data points). Models are optimized using the AdamW \cite{loshchilov2017decoupled} variant of gradient descent on the training set.
For each system, we optimized at least 500 model configurations (obtained via directed and random hyperparameter searches) for 10,000 gradient descent updates. We select the best performing models on the open-loop MSE for the development set, and report performance results on the held-out test set. The hyperparameter space for the model configuration search ranged over i) block structure, ii) linear map prior, iii) nonlinear map type, iv) activation function type,
v) weights on the multi-objective loss function, vi)
gradient descent step size (from \num{3e-5} to 0.01) and $N$-step prediction horizon, vii) estimator input history window, and
viii) neural network width and depth.

\subsection{Results and Discussion}
\begin{figure*}
        \centering
        \includegraphics[width=0.3\textwidth]{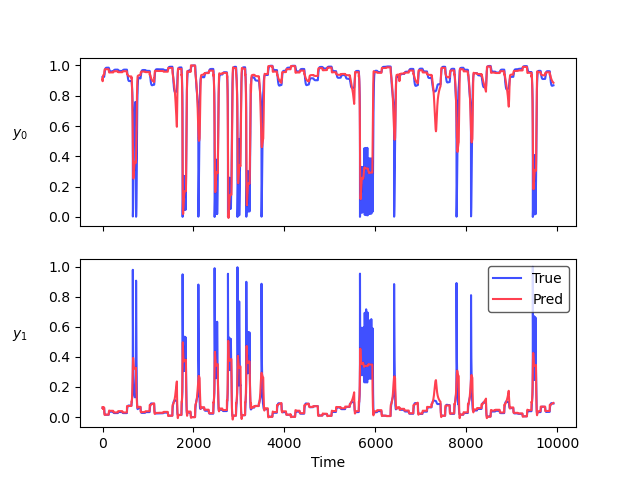}
        \includegraphics[width=0.3\textwidth]{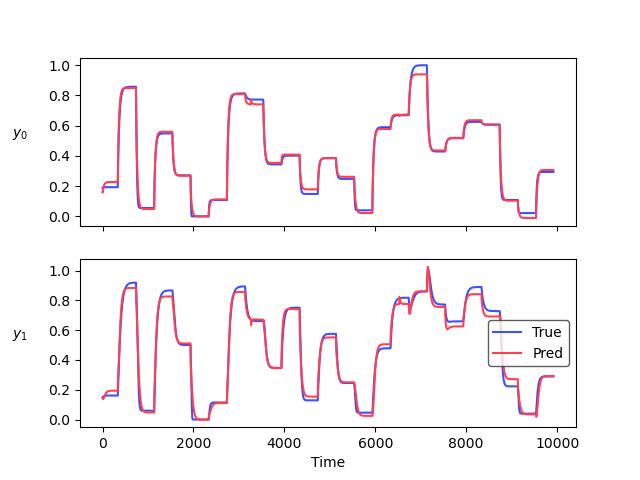}
        \includegraphics[width=0.3\textwidth]{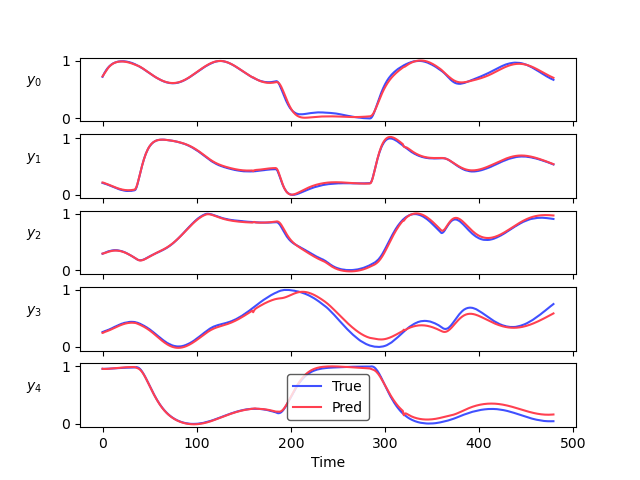}\\
        \includegraphics[width=0.3\textwidth]{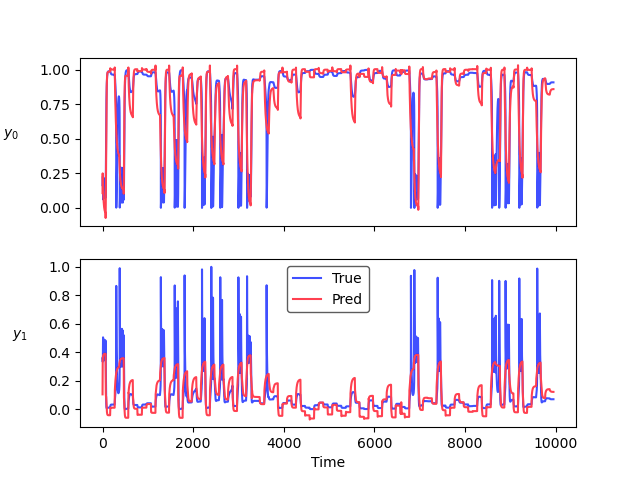}
        \includegraphics[width=0.3\textwidth]{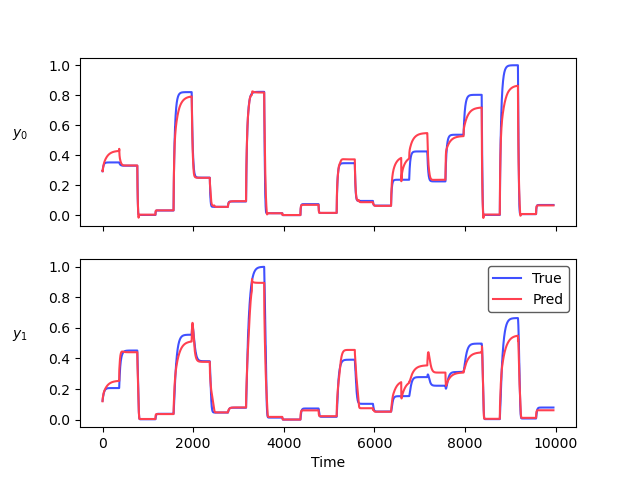}
        \includegraphics[width=0.3\textwidth]{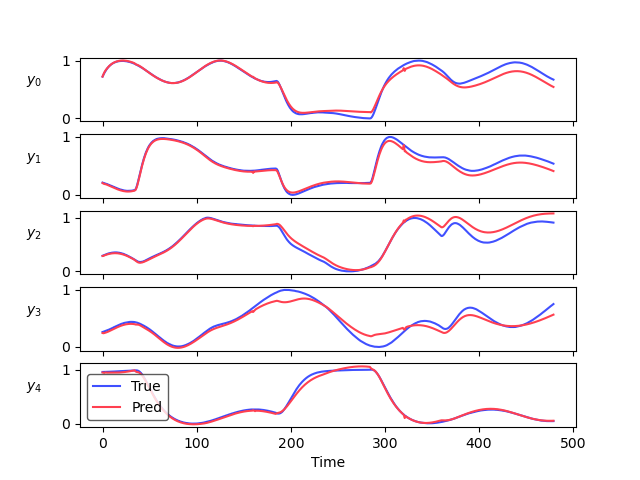}
        \caption{Open-loop traces of the best performing model configurations with constraints (top row) compared to unstructured models without constraints (bottom row) for the CSTR (left), Two Tank (center), and Aerodynamic Body (right) systems.}
        \label{fig:best_blocknlin_open_loop}
        % \vspace{-10pt}
    \end{figure*}

Table \ref{tab:best_results_by_dev} shows a comparison of our best performing systems models versus models trained with a single $N$-step objective loss~\eqref{eq:loss:y} without constraints penalties with an unstructured neural network model. Table \ref{tab:hyperparams} shows configuration details for the best performing models. The CSTR system shows an almost 200\% improvement on open-loop MSE using a block nonlinear structured model and the constraint enforcing multi-objective loss~\eqref{eq:loss:multi}. The Aerodynamic Body system shows greater than 200\% improvement using Spectral regularization for the linear maps of the neural network blocks, a learnable BLU activation function, and penalty constraints via multi-objective loss~\eqref{eq:loss:multi}. The Two Tank system shows an order of magnitude improvement using the neural Hammerstein model with the multi-objective loss~\eqref{eq:loss:multi}. For this system, the structured model generalizes well to the open loop evaluation task from the $N$-step training objective whereas the unstructured model does not. Figure \ref{fig:best_blocknlin_open_loop} compares the traces of these learned models against the ground truth demonstrating that these quantitative gains translate to meaningful improvement in performance. The Two Tank and Aerodynamic Body systems are well-fit on the development set for both the constrained/structured and unconstrained/unstructured models. However, the unconstrained/unstructured models (second row) drift from the true trajectory on the unobserved test set. For the unstable highly nonlinear CSTR system, neither model reaches all the extreme values but the unconstrained/unstructured model has an evident bias for even the non-extreme states across the development and test sets. 

\begin{figure}
        \centering
        \includegraphics[width=\columnwidth]{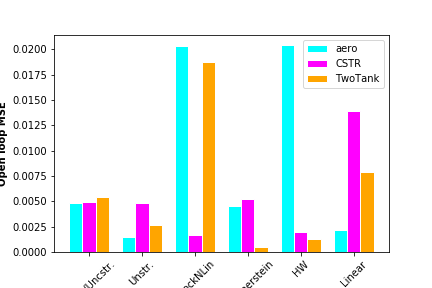}
        \caption{Performance comparison across block structured models.}
        \label{fig:blockcompare}
            \vspace{-10pt}
    \end{figure}

As may be expected, shown in Figure \ref{fig:blockcompare}, particular systems favor particular block structured configurations, with the CSTR model preferring nonlinear maps for all block components, the Two Tank model preferring a Hammerstein model with a nonlinear map only on the input, and the Aerodynamic Body model preferring the model without block structure which allows for more direct nonlinear interaction between the input and state variables.  
    
\begin{table}[t]
\caption{Test set MSE of best-observed structured constrained and unstructured unconstrained models.} 
        \resizebox{\columnwidth}{!}{
        \centering
        \renewcommand{\arraystretch}{1.2}
        \begin{tabular}{llcrr}
        \toprule
          System                    & Structure                          & Constr.           & $N$-step               & Open-loop \\ \midrule
          \multirow{2}{*}{CSTR}     & Block-nonlinear    & Y   & 0.00126   & 0.00171  \\
            & Unstructured\cellcolor{Gray} & N\cellcolor{Gray} & 0.00477\cellcolor{Gray} & 0.00496\cellcolor{Gray} \\\midrule
          \multirow{2}{*}{TwoTank}  & Hammerstein             & Y    & 0.00029     & 0.00050  \\
  & Unstructured\cellcolor{Gray} & N\cellcolor{Gray} & 0.00047\cellcolor{Gray} & 0.00542\cellcolor{Gray}\\\midrule
          \multirow{2}{*}{Aero}    & Unstructured  & Y      & 0.00144     & 0.00153  \\
           & Unstructured\cellcolor{Gray} & N\cellcolor{Gray} & 0.00481\cellcolor{Gray} & 0.00484\cellcolor{Gray} \\
        \bottomrule
        \end{tabular}
        }
        %Results for systems marked with an asterisk are averaged over 5 simulations of the system.}
        \label{tab:best_results_by_dev}
            \vspace{-10pt}
    \end{table}

    \begin{table*}[th!]
        \centering
         \caption{Configurations of best-observed models for each system.}
        \renewcommand{\arraystretch}{1.2}
        \resizebox{\textwidth}{!}{
        \begin{tabular}{lllllllllllllllllll}\toprule
        System      & Structure              & State Est. & LM     & NLM  & Layers & Nodes & Act. & LR         & $Q_{dx}$ & $Q^{\text{con}}_y$ & $Q^{\text{con}}_{\boldsymbol{f}_u}$ & $Q_y$ & $\lambda_\text{min}$ & $\lambda_\text{max}$ & $N$ & $N_p$ \\ \midrule
         CSTR       & Block-nonlinear & RNN        & Linear        & rMLP & 6      & 60    & GELU & \num{1e-4} & 0.2      & 0.2                & 0.1                                 & 0.5   & ---                  & ---                  & 64  & 1 \\
         TwoTank    & Hammerstein     & rMLP       & Linear        & rMLP & 3      & 40    & GELU & \num{3e-4} & 0.3      & 0.3                & 0.1                                 & 0.5   & ---                  & ---                  & 64  & 4 \\
         Aero       & Unstructured    & RNN        & Spectral      & rMLP & 2      & 25    & BLU  & 0.01       & 0.2      & 0.2                & ---                                 & 0.5   & 0.4                  & 0.7                  & 16  & 2 \\
         %Aero       & Unstructured    & MLP        & PF     & RNN  & 5      & 225   & SoftExp & \num{3e-3} & 0.3      & 0.2                & 0.3                                 & 2.0   & \num{1e-5}           & 1.3                  & 16  & 1 \\
            \bottomrule
        \end{tabular}
        }
        \label{tab:hyperparams}
            \vspace{-10pt}
    \end{table*}
    
\subsection{Ablation Study}
The above results suggest that constraint-enforcing objective terms and structural priors can have a dramatic effect on the ultimate model fit. In this section we perform studies to assess which modeling priors have greatest effect for each respective system.    
For each system, we conduct an ablation study in order to assess the effectiveness of constraint penalty terms and linear map parameterizations for modeling. For each system's best-observed architecture, we train 10 randomly initialized models with all the recorded constraint values, no constraints, and ablating each constraint in turn. Specifically, we train models setting $Q^{con}_{y}=0$, $Q_{dx}=0$, $Q^{con}_{f_u}=0$, and with no inductive prior on the linear map. As an extreme case, we also show results for models trained without block structure or constraints penalties. 

\begin{figure}
    \centering
    \includegraphics[width=\columnwidth]{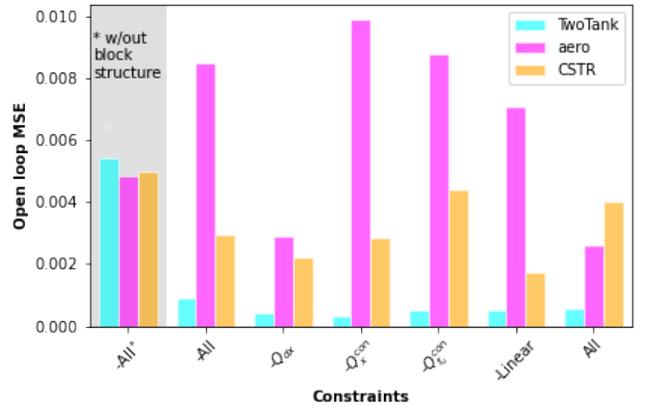}
    \caption{Performance comparison across models with ablated constraints.}
    \label{fig:ablation}
        \vspace{-10pt}
\end{figure}

Figure \ref{fig:ablation} shows the test set open-loop MSE on the best performing models for each system given the different ablation conditions.  
Fully constrained and block-structured nonlinear models (All) perform better than unstructured unconstrained models (-All*).
However, \emph{different} combinations of constraints may improve performance depending on the system under consideration. For the Aerodynamic Body system, removing any one constraint degrades performance similar to removing all constraints with the exception of smoothing ($Q_{dx}$) which seems to have less effect. 
For the Two Tank system removing most individual constraints has little effect but removing the smoothing constraint can actually improve performance which aligns with the discontinuities in the underlying system. Performance on Two Tank appears to rely heavily on the Hammerstein block structure which results in a dramatic drop in performance when removed. The CSTR model performs worse when removing the input influence constraint ($Q^{con}_{f_u}$) or the block-nonlinear model structure, but can benefit when removing other constraints. Similar to the Two Tank model, the near discontinuities in the underlying system behavior make the smoothing term an ill-posed regularization constraint for this particular system. 
% We reserve more detailed analysis of the effect of constraints on optimization and performance to future work. 
% TODO: to add some short speculations?
% \begin{table*}
%     \centering
%     \label{tab:bestablation}
% \begin{tabular}{lrrrrrrr}
% \toprule
% System & Unstr./Unconstr. & Unconstr. &  -$Q_{dx}$ & -$Q^{con}_{x}$ & -$Q^{con}_{f_u}$ & -linear & All constr.\\
% \midrule
% TwoTank &    0.00542 &   0.00090 & 0.00040 & 0.00032 &   0.00050 &    0.00050 &  0.00058 \\
% aero    &    0.00484 &   0.00846 & 0.00291 & 0.00988 &   0.00877 &    0.00708 &  0.00260 \\
% CSTR    &    0.00496 &   0.00295 & 0.00221 & 0.00285 &   0.00440 &    0.00171 &  0.00399 \\
% \bottomrule
% \end{tabular}
% \caption{Caption}
% \end{table*}
% ANALYSIS:
%  1, performance of best nonlinear models (both structured and unstructured)
%  2, ablation 2: comparison without constraints
%  3, weight eigenvalues analysis - link with dissipativness and energy conversion of the system
%  4, optional: plot learnable activation functions - each block with different colouring

    \subsection{Eigenvalue Analysis} 
    \begin{figure*}[t!]
        \centering
        \includegraphics[width=0.23\textwidth]{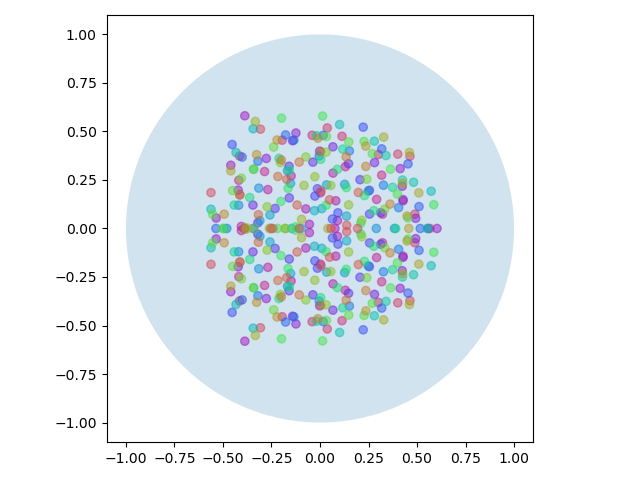}
        \includegraphics[width=0.23\textwidth]{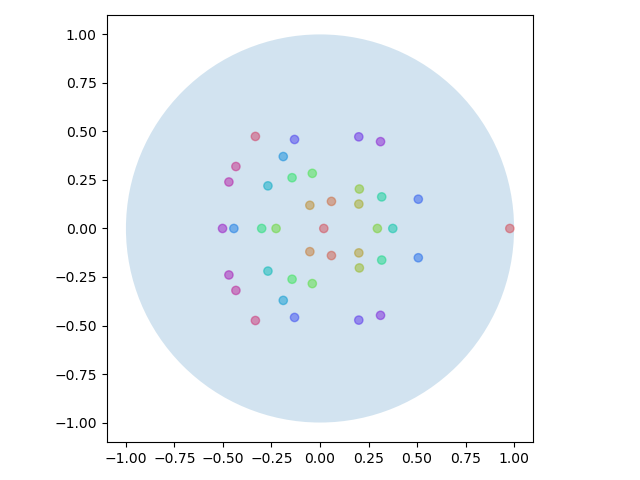}
        \includegraphics[width=0.23\textwidth]{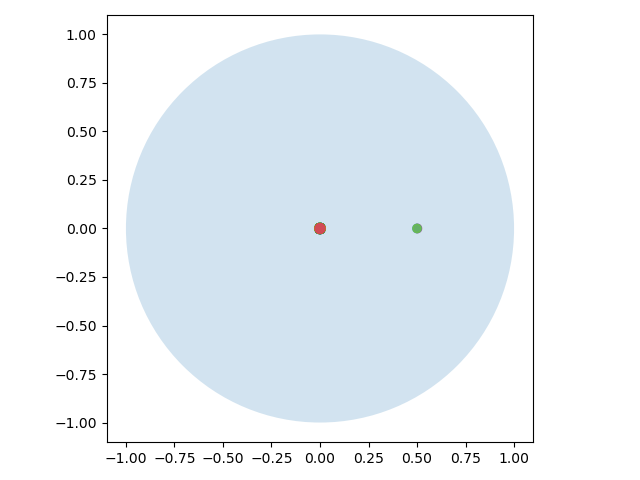}\\
        \includegraphics[width=0.23\textwidth]{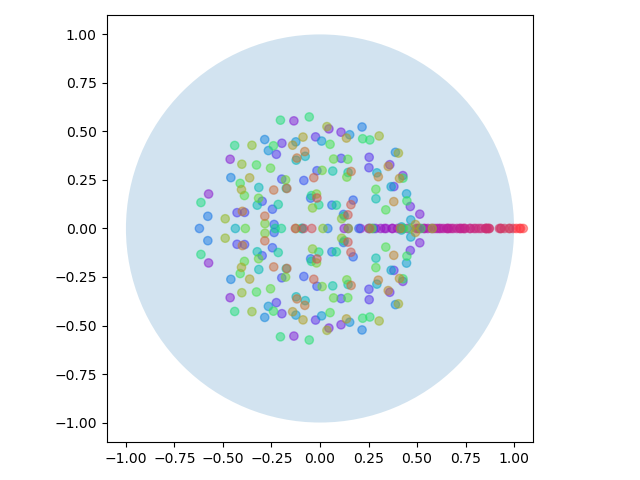}
        \includegraphics[width=0.23\textwidth]{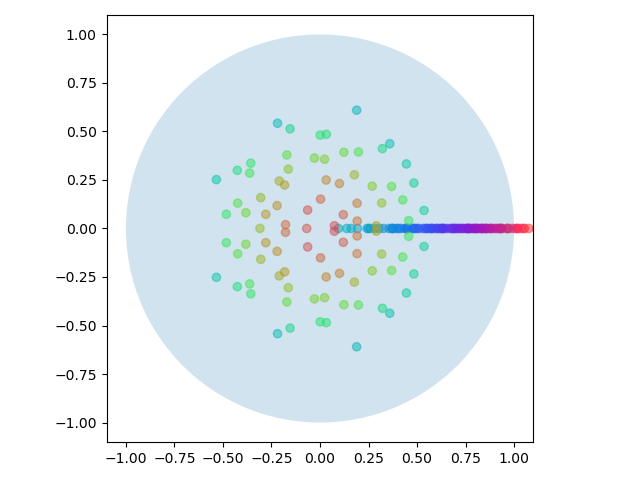}
        \includegraphics[width=0.23\textwidth]{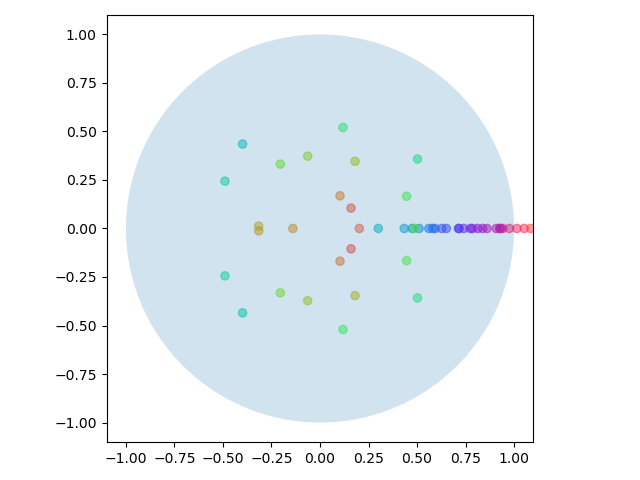}\\
        \includegraphics[width=0.23\textwidth]{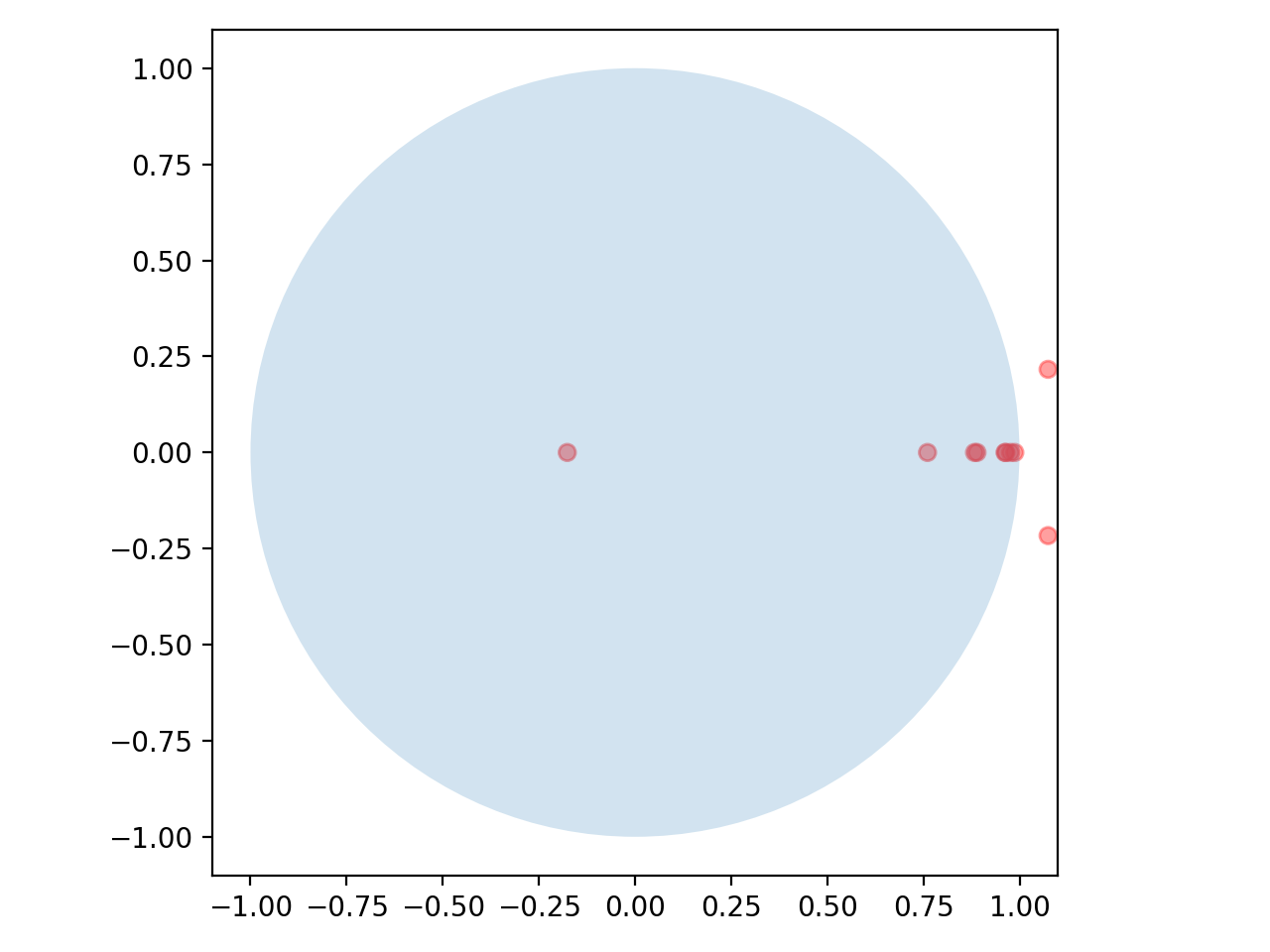}
        \includegraphics[width=0.23\textwidth]{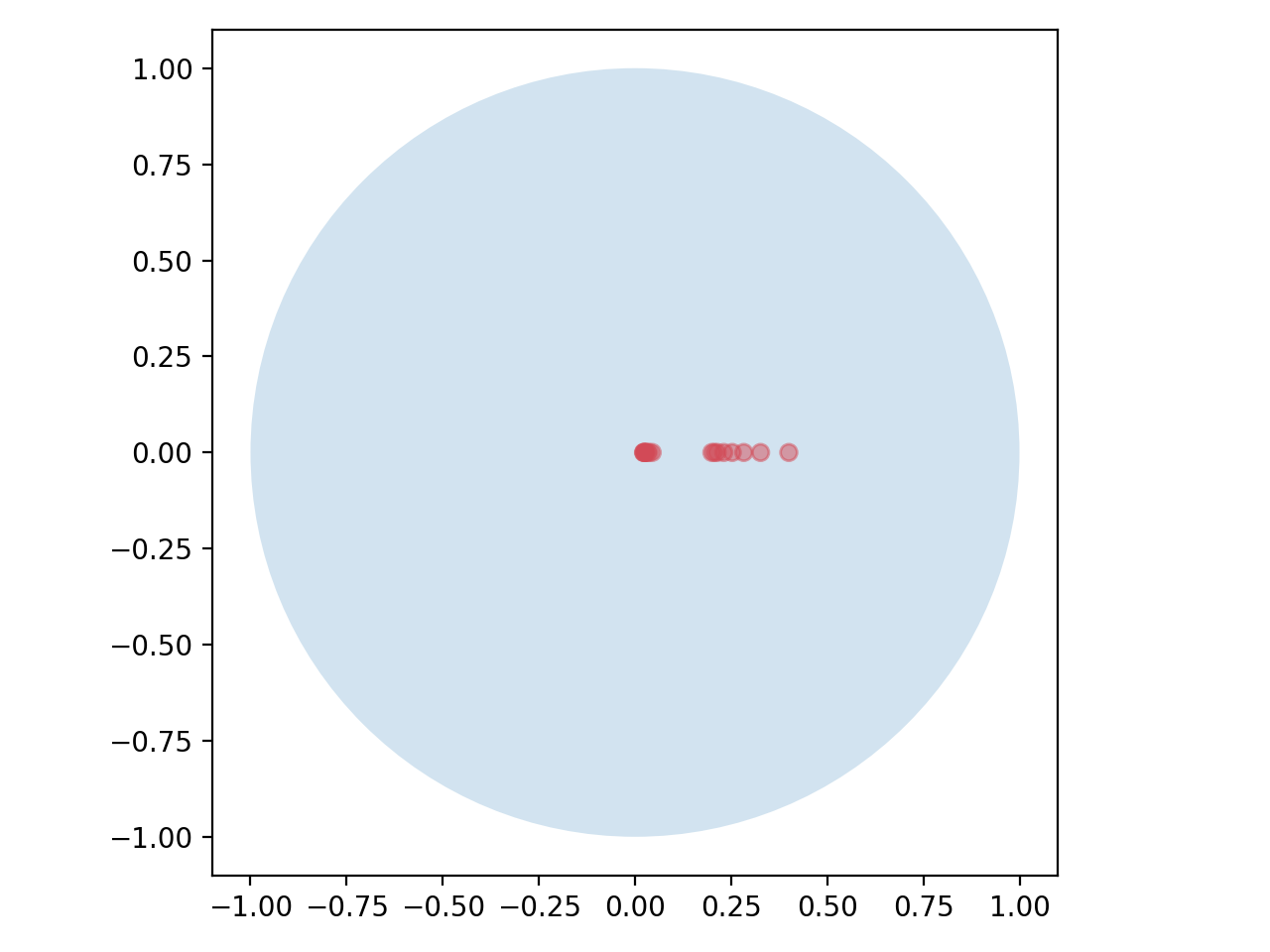}
        \includegraphics[width=0.23\textwidth]{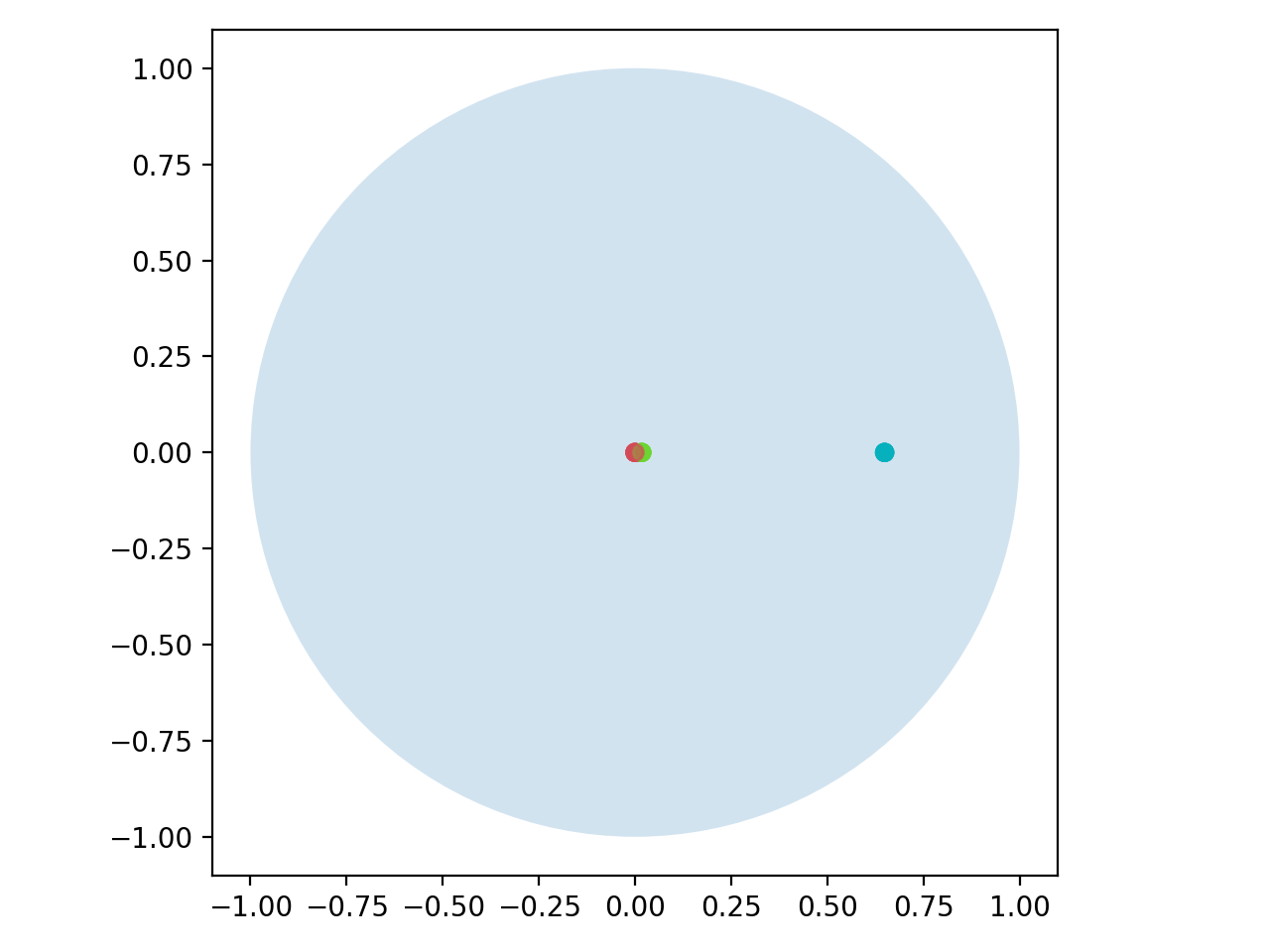}
        \caption{Eigenvalue plots in the complex plane for the state transition map weights of the best performing models with constraints (first row) compared to unstructured and unconstrained (second row) for the CSTR (left), Two Tank (center), and Aerodynamic Body (right) systems. 
        The third row plots eigenvalues of ground truth ODEs linearized across a range of initial conditions.
        Blue circles represent stable regions.}
        \label{fig:eigenvalues}
            \vspace{-10pt}
    \end{figure*}
    Fig. \ref{fig:eigenvalues} shows concatenated eigenvalue plots of the state transition maps of the best performing structured and constrained (top row) vs unstructured and unconstrained models  (second row) for the CSTR, Two Tank, and Aerodynamic Body systems, respectively. We display only the eigenvalues of the linear maps for $\mathbf{f}_x$ and  $\mathbf{f}_{xu}$, respectively. Hence the damping effect of the used stable activation functions on the eigenvalue space is omitted in this analysis. Blue circles represent the stable regions of the complex plane. 
    An interesting observation is that the best performing constrained and structured models (first row) always yield stable state transition weights. 
    By intuition, the composition of stable functions yields the global stability of the learned system dynamics. 
    A formal proof of this proposition will be subject to future work.
    On the other hand, the unstructured and unconstrained models (second row) sometimes learn weights with unstable eigenvalues, therefore, lose the global stability guarantee. 
    
    The eigenvalues of the linearized ground truth ODEs (bottom row) show the tendencies of each system over a range of steady-state conditions. The unstable oscillatory modes of CSTR are apparent in the original linearized system (bottom left) and explain why the stable eigenvalue spectrum of the learned block-structured model (top left) is not entirely capturing the oscillatory modes, as shown in Fig.~\ref{fig:best_blocknlin_open_loop}.
    The small eigenvalues are expected for the Two Tank system (bottom center) because the system's energy drops rapidly without inlet flow.
    In contrast, the learned Hammerstein model (top center) has a dominant eigenvalue close to 1, which refers to the energy conversion property of the learned latent dynamics. However, thanks to the projection of the latent dynamics to the output space dynamics via linear $\mathbf{f}_y$ map, the model is still accurately representing the actual system dynamics. 
    The Aerodynamic Body model is linearized around its steady trimmed conditions and the resulting eigenvalues (bottom right) are very similar to the ones of the learned constrained block-structured model (top right).

\section{Conclusion}
We proposed constrained block nonlinear neural state space models for the identification of nonlinear dynamical systems, and a host of options for constructing such models. Analyses of three dynamical systems have exhibited the efficacy of this modeling approach by superior open-loop performance. In addition to learning accurate systems models from small measurement sets, our constrained and block-structured neural models can yield meaningful performance gains over their unconstrained and unstructured counterparts. Different model structures performed best for each test system, emphasizing the utility of our adaptable, modular deep system identification methodology. %Finally, our ablation studies show objective constraints and structural priors can be beneficial, and the best choices of constraints and priors are highly system dependent. 
Our analysis of the learned system dynamics suggests that suitable constraints and structural priors can be chosen for effective data-driven modeling given some high-level intuition of the underlying system and lacking any specific knowledge of the governing equations or their functional forms. 
% We are motivated by these promising results and believe further analysis of constraints and architectures would be useful. 
For future work, extended analysis of individual constraints and neural architecture choices on a broader set of non-autonomous and autonomous systems can provide further guidance on effective application of our presented methods for engineering practice. 
%We are currently conducting a  larger  case  study  exploring ..
%In future work we will explore a broader set of non-autonomous systems as well as autonomous systems.
%to demonstrate accurate system identification on an extended range of nonlinear systems.
% Further work exploring individual constraints effects is proposed

\bibliography{bib}

\begin{thebibliography}{10}

\bibitem{khalil_2002}
Hassan~K Khalil.
\newblock {\em {Nonlinear systems; 3rd ed.}}
\newblock Prentice-Hall, Upper Saddle River, NJ, 2002.

\bibitem{ioann_thesis}
Yani~Andrew Ioannou.
\newblock {\em Structural Priors in Deep Neural Networks}.
\newblock PhD thesis, University of Cambridge, 2018.

\bibitem{schoukens2019}
J.~{Schoukens} and L.~{Ljung}.
\newblock Nonlinear system identification: A user-oriented road map.
\newblock {\em IEEE Control Systems Magazine}, 39(6):28--99, 2019.

\bibitem{bai_introduction_2010}
Er-Wei Bai and Fouad Giri.
\newblock Introduction to {Block}-oriented {Nonlinear} {Systems}.
\newblock In Fouad Giri and Er-Wei Bai, editors, {\em Block-oriented
  {Nonlinear} {System} {Identification}}, Lecture {Notes} in {Control} and
  {Information} {Sciences}, pages 3--11. Springer, London, 2010.

\bibitem{brockett_volterra_1976}
Roger~W. Brockett.
\newblock Volterra series and geometric control theory.
\newblock {\em Automatica}, 12(2):167--176, March 1976.

\bibitem{oh2011application}
Hyun-Joo Oh and Biswajeet Pradhan.
\newblock Application of a neuro-fuzzy model to landslide-susceptibility
  mapping for shallow landslides in a tropical hilly area.
\newblock {\em Computers \& Geosciences}, 37(9):1264--1276, 2011.

\bibitem{chiras2001}
N.~{Chiras}, C.~{Evans}, and D.~{Rees}.
\newblock Nonlinear gas turbine modeling using narmax structures.
\newblock {\em IEEE Transactions on Instrumentation and Measurement},
  50(4):893--898, 2001.

\bibitem{kibangou2006}
A.~Y. {Kibangou} and G.~{Favier}.
\newblock Wiener-hammerstein systems modeling using diagonal volterra kernels
  coefficients.
\newblock {\em IEEE Signal Processing Letters}, 13(6):381--384, 2006.

\bibitem{gilabert2005}
P.~{Gilabert}, G.~{Montoro}, and E.~{Bertran}.
\newblock On the wiener and hammerstein models for power amplifier
  predistortion.
\newblock In {\em 2005 Asia-Pacific Microwave Conference Proceedings},
  volume~2, pages 4 pp.--, 2005.

\bibitem{Bai2003}
{Er-Wei Bai}.
\newblock Frequency domain identification of hammerstein models.
\newblock {\em IEEE Transactions on Automatic Control}, 48(4):530--542, 2003.

\bibitem{dempsy2004}
E.~J. {Dempsey} and D.~T. {Westwick}.
\newblock Identification of hammerstein models with cubic spline
  nonlinearities.
\newblock {\em IEEE Transactions on Biomedical Engineering}, 51(2):237--245,
  2004.

\bibitem{kwak1998}
{Byung-Jae Kwak}, A.~E. {Yagle}, and J.~A. {Levitt}.
\newblock Nonlinear system identification of hydraulic actuator. friction
  dynamics using a hammerstein model.
\newblock In {\em Proceedings of the 1998 IEEE International Conference on
  Acoustics, Speech and Signal Processing, ICASSP '98 (Cat. No.98CH36181)},
  volume~4, pages 1933--1936 vol.4, 1998.

\bibitem{pearson_gray-box_2000}
Ronald~K. Pearson and Martin Pottmann.
\newblock Gray-box identification of block-oriented nonlinear models.
\newblock {\em Journal of Process Control}, 10(4):301--315, August 2000.

\bibitem{chen1990}
S.~Chen, S.~A. Billings, and P.~M. Grant.
\newblock Non-linear system identification using neural networks.
\newblock {\em International Journal of Control}, 51(6):1191--1214, 1990.

\bibitem{chow1998recurrent}
Tommy~WS Chow and Yong Fang.
\newblock A recurrent neural-network-based real-time learning control strategy
  applying to nonlinear systems with unknown dynamics.
\newblock {\em IEEE transactions on industrial electronics}, 45(1):151--161,
  1998.

\bibitem{al2008nonlinear}
RK~Al~Seyab and Yi~Cao.
\newblock Nonlinear system identification for predictive control using
  continuous time recurrent neural networks and automatic differentiation.
\newblock {\em Journal of Process Control}, 18(6):568--581, 2008.

\bibitem{yu2004nonlinear}
Wen Yu.
\newblock Nonlinear system identification using discrete-time recurrent neural
  networks with stable learning algorithms.
\newblock {\em Information sciences}, 158:131--147, 2004.

\bibitem{de2007nonlinear}
Jos{\'e} de~Jes{\'u}s~Rubio and Wen Yu.
\newblock Nonlinear system identification with recurrent neural networks and
  dead-zone kalman filter algorithm.
\newblock {\em Neurocomputing}, 70(13-15):2460--2466, 2007.

\bibitem{LSTM_SysID2017}
{Yu Wang}.
\newblock A new concept using lstm neural networks for dynamic system
  identification.
\newblock In {\em 2017 American Control Conference (ACC)}, pages 5324--5329,
  2017.

\bibitem{woo_dynamic_2018}
Joohyun Woo, Jongyoung Park, Chanwoo Yu, and Nakwan Kim.
\newblock Dynamic model identification of unmanned surface vehicles using deep
  learning network.
\newblock {\em Applied Ocean Research}, 78:123--133, September 2018.

\bibitem{wilson_neural_2009}
Robert Wilson and Leif Finkel.
\newblock A {Neural} {Implementation} of the {Kalman} {Filter}.
\newblock In Y.~Bengio, D.~Schuurmans, J.~D. Lafferty, C.~K.~I. Williams, and
  A.~Culotta, editors, {\em Advances in {Neural} {Information} {Processing}
  {Systems} 22}, pages 2062--2070. Curran Associates, Inc., 2009.

\bibitem{coskun2017long}
Huseyin Coskun, Felix Achilles, Robert DiPietro, Nassir Navab, and Federico
  Tombari.
\newblock Long short-term memory kalman filters: Recurrent neural estimators
  for pose regularization.
\newblock In {\em 2017 IEEE International Conference on Computer Vision
  (ICCV)}, pages 5525--5533. IEEE, 2017.

\bibitem{krishnan2015deep}
Rahul~G Krishnan, Uri Shalit, and David Sontag.
\newblock Deep kalman filters.
\newblock {\em stat}, 1050:25, 2015.

\bibitem{yeung2019learning}
Enoch Yeung, Soumya Kundu, and Nathan Hodas.
\newblock Learning deep neural network representations for koopman operators of
  nonlinear dynamical systems.
\newblock In {\em 2019 American Control Conference (ACC)}, pages 4832--4839.
  IEEE, 2019.

\bibitem{chen2018neural}
Ricky~TQ Chen, Yulia Rubanova, Jesse Bettencourt, and David~K Duvenaud.
\newblock Neural ordinary differential equations.
\newblock In {\em Advances in neural information processing systems}, pages
  6571--6583, 2018.

\bibitem{champion2019data}
Kathleen Champion, Bethany Lusch, J~Nathan Kutz, and Steven~L Brunton.
\newblock Data-driven discovery of coordinates and governing equations.
\newblock {\em Proceedings of the National Academy of Sciences},
  116(45):22445--22451, 2019.

\bibitem{Raissi_PINN_p1_2017}
Maziar Raissi, Paris Perdikaris, and George~E. Karniadakis.
\newblock Physics informed deep learning (part {I):} data-driven solutions of
  nonlinear partial differential equations.
\newblock {\em CoRR}, abs/1711.10561, 2017.

\bibitem{raissi_multistep_2018}
Maziar Raissi, Paris Perdikaris, and George~Em Karniadakis.
\newblock Multistep {Neural} {Networks} for {Data}-driven {Discovery} of
  {Nonlinear} {Dynamical} {Systems}.
\newblock {\em arXiv:1801.01236 [nlin, physics:physics, stat]}, January 2018.
\newblock arXiv: 1801.01236.

\bibitem{HeZRS15}
Kaiming He, Xiangyu Zhang, Shaoqing Ren, and Jian Sun.
\newblock Deep residual learning for image recognition.
\newblock {\em CoRR}, abs/1512.03385, 2015.

\bibitem{NIPS2018_7892_NeuralODEs}
Tian~Qi Chen, Yulia Rubanova, Jesse Bettencourt, and David~K Duvenaud.
\newblock Neural ordinary differential equations.
\newblock In S.~Bengio, H.~Wallach, H.~Larochelle, K.~Grauman, N.~Cesa-Bianchi,
  and R.~Garnett, editors, {\em Advances in Neural Information Processing
  Systems 31}, pages 6571--6583. Curran Associates, Inc., 2018.

\bibitem{DeepXDE2019}
Lu~Lu, Xuhui Meng, Zhiping Mao, and George~E. Karniadakis.
\newblock {DeepXDE:} {A} deep learning library for solving differential
  equations.
\newblock {\em CoRR}, abs/1907.04502, 2019.

\bibitem{IMEXnet2019}
Eldad Haber, Keegan Lensink, Eran Treister, and Lars Ruthotto.
\newblock {IMEXnet:} {A} forward stable deep neural network.
\newblock {\em CoRR}, abs/1903.02639, 2019.

\bibitem{HaberR17}
Eldad Haber and Lars Ruthotto.
\newblock Stable architectures for deep neural networks.
\newblock {\em CoRR}, abs/1705.03341, 2017.

\bibitem{krishnan2016structured}
Rahul~G. Krishnan, Uri Shalit, and David Sontag.
\newblock Structured inference networks for nonlinear state space models, 2016.

\bibitem{LatentDynamics2018}
Danijar Hafner, Timothy~P. Lillicrap, Ian Fischer, Ruben Villegas, David Ha,
  Honglak Lee, and James Davidson.
\newblock Learning latent dynamics for planning from pixels.
\newblock {\em CoRR}, abs/1811.04551, 2018.

\bibitem{MastiCDC2018}
D.~{Masti} and A.~{Bemporad}.
\newblock Learning nonlinear state-space models using deep autoencoders.
\newblock In {\em 2018 IEEE Conference on Decision and Control (CDC)}, pages
  3862--3867, 2018.

\bibitem{NIPS2018_8004}
Syama~S. Rangapuram, Matthias~W. Seeger, Jan Gasthaus, Lorenzo Stella, Yuyang
  Wang, and Tim Januschowski.
\newblock Deep state space models for time series forecasting.
\newblock In S.~Bengio, H.~Wallach, H.~Larochelle, K.~Grauman, N.~Cesa-Bianchi,
  and R.~Garnett, editors, {\em Advances in Neural Information Processing
  Systems 31}, pages 7785--7794. Curran Associates, Inc., 2018.

\bibitem{ogunmolu_nonlinear_2016}
Olalekan Ogunmolu, Xuejun Gu, Steve Jiang, and Nicholas Gans.
\newblock Nonlinear {Systems} {Identification} {Using} {Deep} {Dynamic}
  {Neural} {Networks}.
\newblock {\em arXiv:1610.01439 [cs]}, October 2016.
\newblock arXiv: 1610.01439.

\bibitem{OgunmoluGJG16}
Olalekan~P. Ogunmolu, Xuejun Gu, Steve~B. Jiang, and Nicholas~R. Gans.
\newblock Nonlinear systems identification using deep dynamic neural networks.
\newblock {\em CoRR}, abs/1610.01439, 2016.

\bibitem{HW_RNN2008}
{Jeen-Shing Wang} and {Yi-Chung Chen}.
\newblock A hammerstein-wiener recurrent neural network with universal
  approximation capability.
\newblock In {\em 2008 IEEE International Conference on Systems, Man and
  Cybernetics}, pages 1832--1837, Oct 2008.

\bibitem{HendrycksG16}
Dan Hendrycks and Kevin Gimpel.
\newblock Bridging nonlinearities and stochastic regularizers with gaussian
  error linear units.
\newblock {\em CoRR}, abs/1606.08415, 2016.

\bibitem{godfrey2019evaluation}
Luke~B Godfrey.
\newblock An evaluation of parametric activation functions for deep learning.
\newblock In {\em 2019 IEEE International Conference on Systems, Man and
  Cybernetics (SMC)}, pages 3006--3011. IEEE, 2019.

\bibitem{he2015delving}
Kaiming He, Xiangyu Zhang, Shaoqing Ren, and Jian Sun.
\newblock Delving deep into rectifiers: Surpassing human-level performance on
  imagenet classification.
\newblock In {\em Proceedings of the IEEE international conference on computer
  vision}, pages 1026--1034, 2015.

\bibitem{zhang2018stabilizing}
Jiong Zhang, Qi~Lei, and Inderjit Dhillon.
\newblock Stabilizing gradients for deep neural networks via efficient svd
  parameterization.
\newblock In {\em International Conference on Machine Learning}, pages
  5806--5814, 2018.

\bibitem{loshchilov2017decoupled}
Ilya Loshchilov and Frank Hutter.
\newblock Decoupled weight decay regularization.
\newblock {\em arXiv preprint arXiv:1711.05101}, 2017.

\bibitem{aero}
{Modeling an Aerodynamic Body}.
\newblock
  \url{https://www.mathworks.com/help/ident/ug/modeling-an-aerodynamic-body.html}.
\newblock Accessed: 2020-09-22.

\end{thebibliography}
\bibliographystyle{unsrt}

% OPTIONAL FUN actions after first draft:
%  1, Potential tuning on CSTR, TwoTank, Aero
%  2, Dubins experiment runs analysis
%  3, Autonomous systems: BasisLinear as nonlinear map, NonlinearExpansion as component, architectures with time delays, different initial conditions
%  4, Moving horizon open loop simulator instead of openloop simulator for performance assessment
% 5, dataset manipulation: training on derivatives instead of normalized values of input/output pairs
% 6, normalize to [-1, 1]
% 7, for dubins model - output equation y = fy(x,u)

% \section*{Appendix}

    % \begin{figure*}[t]
    %     \centering
    %     \includegraphics[width=0.9\textwidth]{figs/eigenvalues_plot.eps}
    %     \caption{Eigenvalue plots of the state transition maps of the best performing models with constraints (top row) compared to unconstrained (bottom row) for the CSTR (left), Two Tank (center), and Aerodynamic (right) systems. Blue circles represent stable regions.}
    %     \label{fig:eigenvalues}
    % \end{figure*}

\end{document}